\documentclass[12pt]{article}
\usepackage{amsmath,amsthm, amssymb,tikz,enumitem}
\usepackage[authoryear,round,longnamesfirst]{natbib}
\usepackage[hidelinks]{hyperref}

\expandafter\def\csname biburl@B13\endcsname{https://arxiv.org/abs/1109.5668}
\expandafter\def\csname biburl@BCK14\endcsname{https://arxiv.org/abs/1205.6721}
\expandafter\def\csname biburl@B16\endcsname{https://arxiv.org/abs/1406.5660}
\expandafter\def\csname biburl@BK18\endcsname{https://arxiv.org/abs/1708.02134}
\expandafter\def\csname biburl@BT19\endcsname{https://arxiv.org/abs/1805.11989}
\expandafter\def\csname biburl@EKMS00\endcsname{https://annals.math.princeton.edu/articles/11101}
\expandafter\def\csname biburl@JM02\endcsname{https://doi.org/10.1016/S0304-4149(01)00142-9}
\expandafter\def\csname biburl@JM04\endcsname{https://arxiv.org/abs/math/0301055}
\expandafter\def\csname biburl@hambly-martin\endcsname{https://arxiv.org/abs/math/0604189}
\expandafter\def\csname biburl@CGK93\endcsname{https://doi.org/10.1214/aoap/1177005277}
\expandafter\def\csname biburl@GK94\endcsname{https://doi.org/10.1214/aoap/1177005201}
\expandafter\def\csname biburl@ganguly2024last\endcsname{https://arxiv.org/abs/2411.08018}
\expandafter\def\csname biburl@LP18\endcsname{https://doi.org/10.1017/9781316104477}

\usepackage{mathtools}

\newcommand{\E}{\mathbb{E}}
\newcommand{\R}{\mathbb{R}}
\newcommand{\Z}{\mathbb{Z}}
\renewcommand{\P}{\mathbb{P}}

\def\l{\left}
\def\r{\right}

\newtheorem{theorem}{Theorem}

\newtheorem{lemma}{Lemma}

\theoremstyle{remark}
\newtheorem{remark}{\bf Remark}

\begin{document}
\title{Sharp integral condition for linear growth of greedy lattice animals}
\author{Jnaneshwar Baslingker \and Daniel Spivak \and B\'alint Vir\'ag}
\date{}
\maketitle

\begin{abstract}
We show that the necessary and sufficient condition for linear growth of greedy lattice animals with nonnegative i.i.d. vertex weights in every dimension \(q\ge2\) is the same criterion for linear growth of directed last-passage times. We also obtain two-sided finite-size estimates, up to dimension-dependent constants,
for the mean maximal weights of directed paths, lattice animals, valid for every nonnegative weights.
These estimates determine the order of mean growth even when growth
is superlinear. Since every up-right path from the origin is a lattice animal, the lower bound follows from the corresponding last-passage bound. For the upper bound, we couple the lattice weights to a Poisson process and represent each animal by a polygonal route of length proportional to its size.
\end{abstract}

\section{Introduction and main results}
Let $q \ge 2$, and let $\{X_\mathbf{v}; \mathbf{v} \in \mathbb{Z}^q\}$ be an i.i.d.\ family of non-negative random variables, with common distribution $F$. For a finite subset $\zeta$ of $\mathbb{Z}^q$, the weight $S(\zeta)$ of $\zeta$ is defined by
\begin{equation}
S(\zeta) = \sum_{\mathbf{v} \in \zeta} X_\mathbf{v}.
\end{equation}
A ``greedy lattice animal" of size $n$ is a connected subset of $\mathbb{Z}^q$ of size $n$ containing the origin, whose weight is maximal among all such sets. Let $N(n)$ be this maximum weight.

This model is studied by \citet{CGK93}, and has applications in statistical physics, queueing theory and percolation. Under the condition that
\begin{equation*}
\mathbb{E}[X_\mathbf{0}^q (\log^+ X_\mathbf{0})^{q+\epsilon}] < \infty \quad \text{for some } \epsilon > 0, 
\end{equation*}
\cite{CGK93} shows that there exists an $N < \infty$ such that
\begin{equation*}
\limsup_{n \to \infty} \frac{N(n)}{n} \le N \quad \text{almost surely.} 
\end{equation*}
\cite{GK94} proceed to show that, under the same condition, one in fact has more strongly that there exists an $N(X) < \infty$ such that
\begin{equation}\label{eq: a.s. limit lattice animal}
\frac{N(n)}{n} \to N \quad \text{almost surely and in } L^1.
\end{equation}
 \cite{JM02} showed the same with the weaker assumption \begin{equation}\label{eq: martin bound}
\int_{0}^{\infty}\l(\P\l(X_{\mathbf{0}}>t\r)\r)^{1/q}dt<\infty.
 \end{equation}

 In the other direction, \cite{CGK93} showed that \(\mathbb E X_{\mathbf 0}^{q}<\infty\) is necessary for linear growth. In dimension two, the work of \cite{ganguly2024last} shows that for directed last-passage percolation, this moment condition is not sufficient, since every directed path is a lattice animal. The present work identifies the necessary and sufficient tail condition in every dimension \(q\ge2\).

\cite{JV26, V25} extends law of large numbers for greedy animals to an ergodic and Poissonian environment. Greedy continuous paths were introduced by
\cite{GM08} as a tool for the study of a continuous model of first-passage percolation. They proved that under an assumption similar to \eqref{eq: martin bound}, the
mass of a greedy continuous path grows linearly. We study such greedy continuous paths first and then transfer the results for continuous models to lattice animals.

Note that $L_n$, the last passage time from $\mathbf{0}\rightarrow\mathbf{n-1}$, is at most $N(q(n-1)+1).$ Let $\nu$ be the measure of $X_{\mathbf{0}}$ and let $N(X)$ be denoted by $\mu_{\Z^q}(\nu).$ We prove the following theorem.
\begin{theorem}\label{thm: bounds for lattice animal}
Let $J_q(\nu)=\l(\int_{0}^{\infty}\nu(x,\infty)^{2/(q+1)}x^{(q-1)/(q+1)}dx\r)$. There exist universal constants $c_q,C_q>0$ such that
    \[
c_qJ_q(\nu)^{\frac{q+1}{2q}}\leq \mu_{\Z^q}(\nu)\leq C_qJ_q(\nu)^{\frac{q+1}{2q}}.
    \]
\end{theorem}
By \citet[Corollary 1.7]{JV26}, the above theorem implies \eqref{eq: a.s. limit lattice animal}.
Note that since the passage time $L_n$, from $\mathbf{0}$ to $\mathbf{n-1}$, is bounded above by $N(q(n-1)+1)$, using \cite[Theorem 1]{BSV26} for the lattice model, the lower bound in Theorem \ref{thm: bounds for lattice animal} follows. We are left to prove only the upper bound. 

\begin{remark}
    We can consider also the greedy lattice path model, treated by \citet{GK94}. Let $M(n)$ be the maximal weight of a self-avoiding
path in $\Z^2$ of length $n$ starting at the origin. Note that if $L(n)$ is the last passage time from $(0,0)$ to $(n-1,n-1)$, then $L(n)\leq M(2n-1)\leq N(2n-1).$ Hence Theorem \ref{thm: bounds for lattice animal} also holds for $M(n)/n$.
\end{remark}

Our results can also be applied to obtain the following finite size bounds for both the last passage times and lattice animal.

\begin{theorem}\label{thm: finite size bounds}
 Let $q\geq 2$ and $X$ be a non-negative random variable. Define 
 \[
 I_{n,q}(X)=n(J_q(X\wedge K_{n,q}))^{\frac{q+1}{2q}}+n^q \int_{K_{n,q}}^{\infty}\P(X>t)dt,
 \]
 where $K_{n,q}=\inf\{t\geq 0: \P(X>t)\leq n^{-q}\}.$
 There exist constants $c,\bar{c}>0$ (dependent on $q$) such that for all $n\geq 1$,

 \[
 cI_{n,q}(X)\leq \E L(n), \E N(n)\leq \bar{c}I_{n,q}(X)
 \]
\end{theorem}
The first term in \(I_{n,q}(X)\) of Theorem \ref{thm: finite size bounds} measures the contribution from weights capped at \(K_{n,q}(X)\), with the collection of weights across multiple scales. The second measures the contribution from rare weights above this cutoff and is comparable to the expected largest excess in a box containing \(n^q\) vertices. Theorem \ref{thm: finite size bounds} shows that all the three models, last passage times, greedy lattice animal, greedy lattice path model have the same mean growth order at finite $n$, in addition to the same integral criterion for linear growth.

To describe Theorem \ref{thm: finite size bounds}, consider two examples of last passage times in $\Z^2$ with $L_n(X),L_n(Y)$ when $X=e^{E/2}(1+E)$ and $Y=e^E/(1+E)^2$ where $E\sim \text{Exp}(1)$ are i.i.d vertex weights. This results in $K_{n,2}(X)\sim n(1+2\log n)$ and $K_{n,2}(Y)\sim n^2/(1+2\log n)^2$. For these two examples we have,
\begin{align*}
    nJ_2(X\wedge K_{n,2}(X))^{3/4}&\asymp n (\log n)^{7/4},\quad \quad n^2\int_{K_{n,2}(X)}^{\infty}\P(X>t)dt \asymp n\log n\\
    nJ_2(Y\wedge K_{n,2}(Y))^{3/4}&\asymp n^2/ (\log n)^{2},\quad \quad n^2\int_{K_{n,2}(Y)}^{\infty}\P(Y>t)dt \asymp n^2/\log n
\end{align*}

For $X,Y$ the dominant terms of $I_{n,2}$ in Theorem \ref{thm: finite size bounds} are the first and second respectively. This gives us $\E L_n(X)\asymp n (\log n)^{7/4}$ and $\E L_n(Y)\asymp n^2 /\log n.$ 

These examples of super-linear growth should be seen in the context of results of \citet{ganguly2024last} which obtained the correct logarithmic
exponent in their lower bound. For critical inverse-square tails,
they proved, with high probability,
\[
L_{\Z^2}(\nu,n)\geq
c\,\frac{n(\log n)^{3/4}}{\log\log n}.
\]
\cite[Theorem 5]{BSV26} showed that for critical heavy tailed weights with $\P(X>t)=t^{-2},$ the growth is of order $n(\log n)^{3/4}.$ 

The approach to prove Theorem \ref{thm: bounds for lattice animal} is similar to \cite{BSV26}. To prove the upper bound in Theorem \ref{thm: bounds for lattice animal}, we prove the corresponding upper bound for a continuous version of the lattice animal that we define below. Then we transfer the results to the lattice animal. We stick to $q=2$ case. The arguments can be extended to $q\geq 3$. 

\section{Continuous version of the lattice animal}

Let $\Pi$ be a Poisson point process on $\R^2\times (0,\infty)$ with intensity \(dxdy\,\eta(dw)\), where \(\eta([t,\infty))<\infty\) for every \(t>0\), is a locally finite positive Borel measure. A route \(\sigma\) is an ordered finite list of distinct points, joined by straight segments, starting and ending at the origin. Denote \(\ell(\sigma)\) for its Euclidean length and \(W(\sigma)\) for the total weight collected by such a route. Define
\[
N_P(t, \eta)=\sup_{\ell(\sigma)\le t}W(\sigma).
\]
We first prove the following result and then transfer the bounds to the lattice animal.
\begin{theorem}\label{thm: cont. lattice animal bound}
    For a universal constant $c>0$,
    \[
    \E N_P(t,\eta)\leq ctJ_2(\eta)^{3/4}.
    \]
\end{theorem}
In order to prove Theorem \ref{thm: cont. lattice animal bound}, we use the following lemma which is proved later. Fix an edge $x$ to $y$ of length $R=|y-x|.$ Let $\mathcal{P}$ be a Poisson point process of intensity $a$ with weight $w$. For an ordered list $\pi$ of $r$ new, distinct points, let $\ell(\pi)$ be the length of path from $x$ to $y$ passing through these ordered points $\pi$, joined by straight lines. Define 
\[
W_{a,w,\theta}(x,y)=\sum_{\pi}
\exp\!\left(\theta\bigl(wr-(\ell(\pi)-R)\bigr)\right).
\]
\begin{lemma}\label{lem: one weight insertion}
 For \(\theta>0\), we have
\[
\mathbb EW_{a,w,\theta}(x,y) \leq 2e^{\theta\gamma R},\quad \text{ if } \gamma:=\left(\frac{4\pi a e^{\theta w}}{\theta^2}\right)^{2/3}\leq 1.
\]

\end{lemma}

Let $\sigma$ be a polygonal path with $h$ distinct points ($h+1$ edges) between $x$ and $y$. A refinement \(\tau\) is obtained by inserting distinct points of \(\mathcal{P}\), in any order, into the edges of \(\sigma\). A refinement preserves the order of the old vertices. In other words, deleting the new points recovers the original ordered route. Write \(N_{\text{new}(\tau)}\) for the total number of new inserted points. Define,
\[
W_{a,w,\theta}(\sigma)=\sum_{\tau\text{ refines }\sigma}
\exp\!\left(\theta\left[
w N_{\text{new}(\tau)}-\bigl(\ell(\tau)-\ell(\sigma)\bigr)
\right]\right)
\]

\begin{lemma}\label{lem: one weight insertion 2}
\[
\mathbb EW_{a,w,\theta}(\sigma)\le 2^{h(\sigma)+1} e^{\theta\gamma\ell(\sigma)},
\qquad
\text{ if }\gamma\le1.
\]
\end{lemma}

\begin{proof}[Proof of Theorem \ref{thm: cont. lattice animal bound}]
For the purpose of this proof, we only consider routes of length at most $t.$
First suppose $\nu=\sum_{j=1}^n a_j\delta_{w_j}$ with $a_j>0$ and
$0<w_{j+1}\le w_j/2$ for $j<n$. Put $T=\sum_{j=1}^n a_j^{2/3}w_j^{4/3}$.
Since multiplying every weight by $s$ multiplies $N_P(t,\nu)$ by $s$ and $T$
by $s^{4/3}$, we normalize the weights so that
\begin{equation}\label{eq:norm}
 (8\pi e)^{2/3}T=1.
\end{equation}

Let $S(\sigma)$, $h(\sigma)$ denote the sum of the weights in the route and the number of points in the route respectively. Let $\mathcal S_i$ consist of routes $\pi$ such that $\ell(\pi)\leq t$ and use only weights $w_j$ with $j\leq i$. Define
\[
 \theta_i=w_i^{-1},\qquad
 c_i=(\log 2)\sum_{j=i+1}^n w_j,\qquad
 L_i(\sigma)=S(\sigma)-\ell(\sigma)+c_i h(\sigma).
\]
Define
\[
 Z_i=\sum_{\sigma\in\mathcal S_i}e^{\theta_i L_i(\sigma)},
 \qquad F_i=\frac{1}{\theta_i}\log Z_i .
\]
The empty route gives $Z_i\ge1$. Let $\mathcal S_0$ contain only the
empty route and $F_0=0$. We also have,
\begin{equation}\label{eq:offset}
 c_n=0,\qquad c_i\le w_i\log2,\qquad
 c_{i-1}-c_i=w_i\log2\quad(2\le i\le n).
\end{equation}

Note that removing all points of weight $w_i$ from a route in $\mathcal S_i$
gives a unique route $\sigma\in\mathcal S_{i-1}$. For a refinement $\pi$ inserting $r$ points
of weight $w_i$,
\[
 L_i(\pi)-L_i(\sigma)
 =(w_i+c_i)r-(\ell(\pi)-\ell(\sigma)).
\]

Then we have,
\[
 Z_i\leq \sum_{\sigma\in\mathcal S_{i-1}}
 e^{\theta_i L_i(\sigma)}W_{a_i,w_i+c_i,\theta_i}(\sigma).
\]

Let $\mathcal F_i$ be generated by the points of weight $j\leq i$ with $\mathcal F_0$ trivial. 
Applying
Lemma \eqref{lem: one weight insertion 2} (length restriction only helps) with $w=w_i+c_i$ gives
\[
 \E[W_{a_i,w_i+c_i,\theta_i}(\sigma)\mid\mathcal F_{i-1}]
 \le 2^{h(\sigma)+1}e^{\theta_i A_i\ell(\sigma)},\qquad
 A_i=\left(4\pi a_i e^{\theta_i(w_i+c_i)}w_i^2\right)^{2/3}.
\]
As \eqref{eq:offset} implies $e^{\theta_i(w_i+c_i)}\le2e$, so
\begin{equation}\label{eq:sum}
 \sum_{i=1}^n A_i\le (8\pi e)^{2/3}T=1.
\end{equation}
In particular the hypothesis $A_i\le1$ holds at each level.
For $i\ge2$,
\eqref{eq:offset} gives
$2^{h(\sigma)}e^{\theta_i L_i(\sigma)}
=e^{\theta_i L_{i-1}(\sigma)}$. Hence
\[
 \E[Z_i\mid\mathcal F_{i-1}]
 \le 2e^{\theta_i A_i t}
 \sum_{\sigma\in\mathcal S_{i-1}}e^{\theta_i L_{i-1}(\sigma)}
 \le 2e^{\theta_i A_i t}Z_{i-1}^{\theta_i/\theta_{i-1}}.
\]
The last inequality uses $\sum_j x_j^p\le(\sum_j x_j)^p$ for
$p=\theta_i/\theta_{i-1}\ge1$. For $i=1$ the sum over old routes
is replaced by $1$. Conditional Jensen's inequality now gives
\[
 \E[F_i\mid\mathcal F_{i-1}]\le F_{i-1}+A_it+w_i\log2,
 \qquad 1\le i\le n .
\]

Taking expectations and summing, we obtain from \eqref{eq:sum}
\[
 \E F_n\le t\sum_{i=1}^n A_i+(\log2)\sum_{i=1}^n w_i
 \le t+2(\log2)w_1 .
\]
Since $c_n=0$, every $\sigma\in\mathcal S_n$ satisfies
$S(\sigma)-\ell(\sigma)\le F_n$. Thus $N_P(t,\nu)\le t+F_n$.
Substituting the normalization \eqref{eq:norm} gives
\begin{equation}\label{eq:finite}
 \E N_P(t,\sum_{j=1}^n a_j\delta_{w_j})\le CtT^{3/4}+2(\log2)w_1 .
\end{equation}

For general $\nu$, let $K=tJ_{2}(\nu)^{3/4}>0$ and finite. For $\nu_K=\nu|_{(0,K]}$ and rounding the weights upwards to dyadic values (largest rounded weight is at most $2K$), for some universal constant $c$ we have $T\leq cJ_2(\nu).$ Applying the finite levels estimate we obtained above and letting the lower cut-off decrease to $0$ we get,
\[
\mathbb E N_P(t,\nu_{\le K})\le ctJ_{2}(\nu)^{3/4}+cK\le 2ctJ_{2}(\nu)^{3/4}.
\]
The expected contribution from points of weight at least $K$ is at most 
\begin{align*}
    \frac{\pi t^2}{4}\int\limits_{[K,\infty)}w\nu(dw)&\leq \frac{\pi t^2}{4}K\nu[K,\infty)+\frac{\pi t^2}{4}\int_K^{\infty}\nu[w,\infty)dw\\
&\leq Ct^2 J_2(\nu)^{3/2}/K.
\end{align*} 
As a result, there exists some universal constant such that
\[
\mathbb E N_P(t,\nu)\le ctJ_{2}(\nu)^{3/4}.
\]
\end{proof}

\begin{proof}[Proof of Lemma \ref{lem: one weight insertion}]
    Let $r$ number of points be inserted in the route from $x$ to $y$. Let $x=p_0,p_1,p_2,..,,p_r,p_{r+1}=y$ be the route and let $z_j=p_j-p_{j-1}$ be the successive steps. Let $e=(y-x)/R$ be the unit vector. If $R=0$, take any unit vector. We have
    \[
    \sum\limits_{j=1}^{r+1}|z_j|=\ell(\pi),\quad \sum\limits_{j=1}^{r+1}e\cdot z_j=e\cdot(y-x)=R.
    \]
    Define 
    \[
K(z)=\exp\{-\theta[|z|-(1-\gamma)e\cdot z]\}.
\]
We have,
\[
\prod_{j=1}^{r+1}K(z_j)
=e^{-\theta[\ell(\pi)-(1-\gamma)R]}
=e^{-\theta(\ell(\pi)-R)}e^{-\theta\gamma R}.
\]
The summand in Lemma \ref{lem: one weight insertion} is exactly
\[
e^{\theta[wr-(\ell(\pi)-R)]}
=e^{\theta\gamma R}e^{\theta wr}\prod_{j=1}^{r+1}K(z_j).
\]
Since $\gamma\leq 1,$ we have $K(z)\leq 1.$ We drop the factor \(K(y-p_r)\). This removes the constraint that the steps must end at \(y\). By the Poisson factorial moment formula \cite[Theorem 4.4 and Corollary 4.10]{LP18}, the expected contribution from routes with exactly \(r\) inserted points is now at most
\[
e^{\theta\gamma R}(ae^{\theta w})^r
\int_{(\mathbb R^2)^r}\prod_{j=1}^{r}K(p_j-p_{j-1})\,dp_1\cdots dp_r
=e^{\theta\gamma R}\bigl(ae^{\theta w}\|K\|_1\bigr)^r.
\]
We now find $\|K\|_1$. Let $b=1-\gamma$.
\[
\|K\|_1
=\int_0^{2\pi}\int_0^\infty s\,e^{-\theta s(1-b\cos\phi)}\,ds\,d\phi
=\frac1{\theta^2}\int_0^{2\pi}\frac{d\phi}{(1-b\cos\phi)^2}
=\frac{2\pi}{\theta^2(1-b^2)^{3/2}}.
\]
This gives us,
\[
ae^{\theta w}\|K\|_1
\le \frac{2\pi ae^{\theta w}}{\theta^2\gamma^{3/2}}
=\frac12.
\]
As a result, we have
\[
\mathbb E W_{a,w,\theta}(x,y)
\le e^{\theta\gamma R}\sum_{r=0}^\infty 2^{-r}
=2e^{\theta\gamma R}.\qedhere
\]
\end{proof}

\begin{proof}[Proof of Lemma \ref{lem: one weight insertion 2}]
    Let $x,v_1,v_2,...,v_{h},y$ be the vertices of the polygonal path. For edge \(j\), denote the ordered list of new points by \(\mathbf z_j=(z_{j,1},\ldots,z_{j,r_j})\). Set \(z_{j,0}=v_{j-1}\) and \(z_{j,r_j+1}=v_j\), and define
\[
\ell_j(\mathbf z_j)
=\sum_{s=0}^{r_j}|z_{j,s+1}-z_{j,s}|,\qquad
f_j(\mathbf z_j)
=\exp\!\left\{\theta\left[wr_j-\bigl(\ell_j(\mathbf z_j)-R_j\bigr)\right]\right\}.
\]
If \(r_j=0\), the edge is unchanged, so \(\ell_j=R_j\) and \(f_j=1\).

Contribution of a refinement $\tau$ to \(W_{a,w,\theta}(\sigma)\) is exactly,
\[
\exp\!\left\{\theta\left[wr-\bigl(\ell(\tau)-\ell(\sigma)\bigr)\right]\right\}
=\prod_{j=1}^{h+1} f_j(\mathbf z_j).
\]
Apply the Poisson factorial moment formula \cite[Theorem 4.4 and Corollary 4.10]{LP18} to get,
\[
\begin{aligned}
a^r\int_{(\mathbb R^2)^r}
   \prod_{j=1}^{h+1}f_j(\mathbf z_j)\,
   d\mathbf z_1\cdots d\mathbf z_{h+1}
&=
\prod_{j=1}^{h+1}
\left[
a^{r_j}\int_{(\mathbb R^2)^{r_j}}
f_j(\mathbf z_j)\,d\mathbf z_j
\right].
\end{aligned}
\]
We have
\[
\mathbb E W_{a,w,\theta}(\sigma)
=
\prod_{j=1}^{h+1}
\left[
\sum_{r_j=0}^{\infty}
a^{r_j}\int_{(\mathbb R^2)^{r_j}}
f_j(\mathbf z_j)\,d\mathbf z_j
\right].
\]
The \(j\)-th factor is exactly \(\mathbb E W_{a,w,\theta}(v_{j-1},v_j)\). Applying Lemma \ref{lem: one weight insertion} to each factor yields,
\[
\mathbb E W_{a,w,\theta}(\sigma)
\le \prod_{j=1}^{h+1}2e^{\theta\gamma R_j}
=2^{h+1}e^{\theta\gamma\ell(\sigma)}.
\]
\end{proof}

\begin{proof}[Proof of Theorem \ref{thm: bounds for lattice animal}]

    Note that since the passage time $L_n$, from $\mathbf{0}$ to $\mathbf{n-1}$, is bounded above by $N_n$, using \cite[Theorem 1]{BSV26} for the lattice model, the lower bound in Theorem \ref{thm: bounds for lattice animal} follows. We are left to prove the upper bound. We prove for the $q=2$ case. The same arguments holds for $q\geq 3.$

Let \(u\) be a median of \(X_0\), and put \(Y_v=X_v\mathbf 1_{\{X_v>u\}}\) for $v\in\Z^2$. Since \(X_v\le u+Y_v\), the weights at most \(u\) contribute at most \(nu\) to any lattice animal.

Consider a Poisson process on \(\mathbb R^2\times(0,\infty)\) with intensity \(dz\,\eta(dw)\), where \(\eta=2\nu|_{(u,\infty)}\). For each \(v\in\mathbb Z^2\), let \(M_v\) be the largest mark in the unit square \(v+[-1/2,1/2)^2\), taking \(M_v=0\) if the square is empty. \(M_v\) stochastically dominates \(Y_v\). By sitewise coupling,
\[
\mathbb E N(n)\le nu+\mathbb E N_M(n),
\]
where \(N_M(n)\) denotes the lattice animal maximum with weights \(M_v\).

Now fix a lattice animal \(A\) of size \(n\) (with weights \(M_v\)). Traverse a spanning tree (in$\Z^2$) of \(A\) from the origin and back, using each of its \(n-1\) edges twice. By looking a the order of first visited vertices, such a walk gives an ordered list of distinct vertices of $\Z^2$. The length of a polygonal route in $\R^2$ (picking Poisson weights) corresponding to that ordered list has length at most $2(n-1)+\sqrt{2}n$. Taking the maximum over \(A\), we obtain the pathwise comparison
\[
N_M(n)\le N_P\bigl((2+\sqrt2)n,\eta\bigr).
\]
As a result, Theorem \ref{thm: cont. lattice animal bound} gives
\[
\mathbb E N(n)\le nu+CnJ_2(\eta)^{3/4}.
\]
Since \(u\le CJ_2(\nu)^{3/4}\), we get the desired bound.
\end{proof}

\begin{remark} [\textbf{Higher dimensions}]
Let \(q\ge2\) and let \(\kappa_q\) be the volume of the Euclidean unit ball in \(\mathbb R^q\). In Lemma \ref{lem: one weight insertion}, take
\[
\gamma=\left(\frac{2q!\kappa_q\,a e^{\theta w}}{\theta^q}\right)^{2/(q+1)}
\le1.
\]
Then the higher dimensional analog of  Lemma \ref{lem: one weight insertion} gives, \(\mathbb E W_{a,w,\theta}(x,y)\le 2e^{\theta\gamma|y-x|}\).  The polygonal bound still has the factor \(2^{h(\sigma)+1}\), since it depends only on the number of edges.

In the upper bound proof, at level \(i\), substituting \(\theta_i=w_i^{-1}\) and \(w=w_i+c_i\) gives
\[
A_i\le C_q a_i^{2/(q+1)}w_i^{2q/(q+1)}.
\]
Thus the same recursion and weight rescaling give, with \(\alpha_q=(q+1)/(2q)\),
\[
\mathbb E N_P(t,\nu)\le
C_qt\left(\sum_i a_i^{2/(q+1)}w_i^{2q/(q+1)}\right)^{\alpha_q}
+C_qw_1.
\]
For general \(\nu\), the same arguments of $q=2$ with \(K=tJ_q(\nu)^{\alpha_q}\) give \(\mathbb E N_P(t,\nu)\le C_qtJ_q(\nu)^{\alpha_q}\). Finally, in the lattice transfer, replace unit squares by unit cubes to get the comparison \(N_M(n)\le N_P((2+\sqrt q)n,\eta)\). The median comparison is unchanged, with \(u\le C_qJ_q(\nu)^{\alpha_q}\).
\end{remark}

\section{Proof of finite size bounds}

We will prove Theorem \ref{thm: finite size bounds} for $q=2$. Similar arguments hold for all $q\geq 3$. Since $ L(n)\leq N(q(n-1)+1)$, in order to prove Theorem \ref{thm: finite size bounds}, we prove the lower bound for $\E L(n)$ and upper bound for $\E N(n).$ We use the notation $L_n(X),N_n(X)$ for the last passage time from $\mathbf{0}\rightarrow\mathbf{n-1}$ and lattice animal of size $n$ respectively, with i.i.d. $X$ vertex weights. We use the following lemma in the proof of the lower bound which is proved later.

\begin{lemma}\label{lem: bounded passage time moment}
Let $Y$ be a bounded random variable with \(0\le Y\le K\) and let \(L_n(Y)\) be the passage time from $(0,0)$ to $(n-1,n-1)$ with i.i.d. \(Y\) vertex weights. Then there exist constants $c_1,c_2>0$ such that, for every \(n\ge1\),
\[
 n\mu_{\text{LPP}}(Y)\le c_1\,\mathbb E L_n(Y)+c_2K,
 \qquad
 \mu_\text{LPP}(Y)=\lim_{r\to\infty}\frac{\mathbb E L_r(Y)}r .
\]
\end{lemma}
\hspace{5mm}
\begin{proof}[Proof of Theorem \ref{thm: finite size bounds}]\mbox{{}}
\newline
\newline
\textbf{Proof of LPP lower bound}

Let \(Y=X\wedge K_n\), where $
K_n=\inf\{t\ge0:\P(X>t)\le n^{-2}\}
$. By the shape lower bound for LPP of \cite[Theorem 1]{BSV26} and Lemma \ref{lem: bounded passage time moment},
\begin{equation*}
    nJ_2(Y)^{3/4}\le Cn\mu_{\mathrm{LPP}}(Y)
\le C\bigl(\mathbb E[L_n(Y)]+K_n\bigr)
\le C\bigl(\mathbb E[L_n(X)]+K_n\bigr).
\end{equation*}

We will show that $K_n\leq C\E L_n(X)$ and $n^2\int_{K_n}^\infty \P(X>t)\,dt
\le C\mathbb \E L_n(X).$ Since the lower bound of Theorem \ref{thm: finite size bounds}, for $q=2$, is
\[
c\l(nJ_2(Y)^{3/4}+n^2\int_{K_n}^\infty \P(X>t)\,dt\r),
\]
we get the desired lower bound.

We first show that $K_n\leq C\E L_n(X)$. Let \(M_n\) be the largest weight in \(\{1,\ldots,n\}^2\). Since \(L_n(X)\ge M_n\) and $P(M_n\ge K_n)\geq C$, we get that $CK_n\leq \E M_n\leq \E L_n(X)$.

We now show $n^2\int_{K_n}^\infty \P(X>t)\,dt
\le C\mathbb \E L_n(X)$. For $t\geq K_n$ we have,
\[
\begin{aligned}
\mathbb P(M_n>t)
=1-(1-\P(X>t))^{n^2}
\ge (1-e^{-1})n^2\P(X>t).
\end{aligned}
\]
Integrating over \(t\ge K_n\), we obtain
\[
(1-e^{-1})n^2\int_{K_n}^\infty \P(X>t)\,dt
\le \mathbb E[(M_n-K_n)_+]
\le \mathbb E[L_n(X)].
\]
\textbf{Proof of lattice animal upper bound}
For any $K>0$, we have
\[
N_n(X)\leq N_n(X\wedge K)+ \sum\limits_{v\in\Z^2:\| v\|_1\leq n-1}(X_v-K)_+.
\]
Taking expectations and applying upper bound from Theorem \ref{thm: bounds for lattice animal} along with super-additivity argument,
\[
\E N_n(X)\leq CnJ_2(X\wedge K)^{3/4}+2n^2\int_{K}^{\infty}\P(X>t)dt.
\]
For $K=K_n$, this gives the desired upper bound.
\end{proof}

\begin{proof}[Proof of Lemma \ref{lem: bounded passage time moment}]
 Let \(B=L_n(Y)\). We first bound the exponential moment of \(B\) using \cite[Theorem 2.1]{BLM00}. For each site \(v\), let \(B^{(v)}\) be the passage time after replacing the weight \(Y_v\) by zero. Fix a maximizing path \(\pi\) in the original environment. Evaluating the same path after the replacement gives
\[
0\le B-B^{(v)}\le Y_v\,\mathbf 1_{\{v\in\pi\}}.
\]
Thus each decrease is at most \(K\), and their sum is at most the weight of \(\pi\), namely \(B\). As a result \(B/K\) is self-bounding. The self-bounding exponential moment inequality of \cite[Theorem 2.1]{BLM00}, applied at \(\lambda=1\), gives 
\begin{equation}\label{eq: bound 1}
    K\log\mathbb E\!\left[\exp(B/K)\right]
\le (e-1)\mathbb E[B].
\end{equation}

Now partition \(\{1,\ldots,mn\}^2\) into \(m^2\) disjoint boxes of side \(n\). For each box \(Q\), let \(B_Q\) be its corner-to-corner passage time. The variables \(B_Q\) are independent copies of \(B\).

Every up-right path from $\mathbf{1}\rightarrow\mathbf{n}$ visits an up-right sequence \(\Gamma\) of boxes. Inside each visited box \(Q\), its path segment has weight at most \(B_Q\). Therefore
\[
L_{mn}(Y)\le \max_{\Gamma}\sum_{Q\in\Gamma}B_Q.
\]
There are \(\binom{2m-2}{m-1}\le 2^{2m-2}\) possible sequences \(\Gamma\), each containing \(2m-1\) distinct boxes. Bounding the exponential of the maximum by the sum of the exponentials gives
\[
\mathbb E\!\left[\exp\!\left(\frac{L_{mn}(Y)}{K}\right)\right]
\le
2^{2m-2}
\left(\mathbb E\!\left[\exp(B/K)\right]\right)^{2m-1}.
\]
Taking logarithms, applying Jensen’s inequality on LHS and \eqref{eq: bound 1} on the right, we obtain
\[
\mathbb E[L_{mn}(Y)]
\le (2m-2)K\log 2+(2m-1)(e-1)\mathbb E[L_n(Y)].
\]
Divide by \(m\) and let \(m\to\infty\). Since \(\mathbb E[L_{mn}(Y)]/m\to n\mu(Y)\), this proves the lemma.
    
\end{proof}

\paragraph{AI assistance.}
The authors derived the exact tail integral condition in a previous work for last passage percolation. The upper bound proof of that work is modified for lattice animals together with ChatGPT Astra 6.1, which was also used to extend the arguments of growth of passage times for critical heavy tailed weights in the said previous work to obtain finite size bounds and it also helped with final editing and error checking.  

\paragraph{Acknowledgments.}
The authors were supported by B.V.'s NSERC Discovery Accelerator grant.

\bibliography{shape}
\bibliographystyle{plainnat-linked}

\par\vspace{2\baselineskip}
\noindent\begin{minipage}{\textwidth}
\normalsize
\noindent
Jnaneshwar Baslingker. Department of Mathematics, University of Toronto,
40 St. George St., Toronto, Ontario, M5S 2E4, Canada.
\href{mailto:j.baslingker@utoronto.ca}{\texttt{j.baslingker@utoronto.ca}}
\par\medskip
\noindent
Daniel Spivak. Department of Mathematics, University of Toronto,
40 St. George St., Toronto, Ontario, M5S 2E4, Canada.
\href{mailto:daniel.spivak@mail.utoronto.ca}{\texttt{daniel.spivak@mail.utoronto.ca}}
\par\medskip
\noindent
B\'alint Vir\'ag. Department of Mathematics, University of Toronto,
40 St. George St., Toronto, Ontario, M5S 2E4, Canada.
\href{mailto:balint@math.toronto.edu}{\texttt{balint@math.toronto.edu}}
\par
\end{minipage}

\end{document}